\documentclass[english]{ccdconf}
\usepackage{amssymb, amsmath}
\usepackage{array}
\usepackage{multirow}
\usepackage{moreverb}
\usepackage{graphicx,epsfig}
\date{}
\usepackage{caption}
\usepackage{algorithm}
\usepackage{algorithmicx}
\usepackage{algpseudocode}
\usepackage{amsmath}
\usepackage{color}
\usepackage{algpseudocode}
\usepackage{float}
\usepackage{amsthm}
\usepackage{amsfonts}
\usepackage{amssymb}
\usepackage{color}
\usepackage{mathrsfs}
\usepackage{amsmath}
\usepackage{cases}
\begin{document}
\pdfoutput=1

\title{Adaptive dynamic programming-based algorithm for infinite-horizon linear quadratic stochastic optimal control problems}
\author{Heng Zhang}
\affiliation{School of Control Science and Engineering,
	Shandong University, Jinan, 250061, China.
        \email{zhangh2828@163.com}}

\maketitle

\begin{abstract}
This paper investigates an infinite-horizon linear quadratic stochastic (LQS) optimal control problem for a class of  continuous-time stochastic systems. By employing the technique of adaptive dynamic programming (ADP), we propose a novel model-free policy iteration (PI) algorithm. Without needing all information of the system coefficient matrices, the proposed PI algorithm iterates by using the data of the input and system state  collected on a fixed time interval. Finally, a numerical example is presented to demonstrate the feasibility of the obtained algorithm.
\end{abstract}

\keywords{Linear quadratic stochastic optimal control, Policy iteration, Adaptive dynamic programming }

\footnotetext{The author acknowledges the financial support from the NSFC under Grant Nos. 11831010, 61925306 and 61821004, and the NSF of Shandong Province under Grant Nos. ZR2019ZD42 and ZR2020ZD24.}
\section{INTRODUCTION}

The linear quadratic stochastic (LQS) optimal control problem, initiated by Wonham \cite{Wonham1968} has been broadly applied in a lot of fields such as engineering. As is known to all, the continuous-time LQS problem in infinite horizon is closely related to the stochastic algebraic Riccati equation (SARE), which is difficult to solve due to its nonlinear structure. With the in-depth study of the LQS optimal control problem, researchers developed some approximation methods to obtain the solution of the SARE. For instance, Ni and Fang \cite{Ni2013} proposed a PI algorithm to solve the SARE iteratively. With the help of the positive operators, a Newton's method was proposed by Damm and  Hinrichsen \cite{Damm2001} to solve the SARE. However, the above methods need all knowledge of the system, i.e., all parameters of the system have to be known beforehand. In fact, the system matrices are difficult to obtain directly in  applications such as engineering and finance. The methods mentioned above will become invalid if the system coefficient matrices are unknown. Thus, it is of great importance to propose a model-free strategy to solve LQS optimal control problems, without using the information of system matrices.
~\\

For the past decade, adaptive dynamic programming (ADP) (Werbos \cite{Werbo1974}) and reinforcement learning (RL) (Sutton and Barto \cite{SuttonBarto1998}) theories have been broadly used to solve optimal control problems with partially model-free or model-free system dynamics. About the development of deterministic system case, see, e.g., Shi and Wang \cite{ShiWang2021}, Pang et al. \cite{PangJM2020}, Kiumarsi et al. \cite{Kiumarsi2017}, Vamvoudakis et al. \cite{Vamvoudakis2017}, Bian and Jiang \cite{Bian2016},   Palanisamy et al. \cite{Palanisamy2015}, Vrabie et al. \cite{Vrabie2009},  Wei et al. \cite{Weizhangdai2009}, Jiang and Jiang \cite{JiangJiang2012}, Mukherjee et al. \cite{Mukherjee2021} and the references therein.
~\\

Regarding to stochastic optimal control problems, Ge et al. \cite{Ge2021} proposed a model-free methodology to get the optimal policy for a kind of mean-field discrete-time stochastic systems by the method of Q-learning. By the technique of ADP, Wang et al. \cite{WangZhangLuo2016} solved a class of discrete-time LQS optimal control problems. Wang et al. \cite{WangZhangLuo2018} developed a model-free Q-learning algorithm to get the optimal control for  discrete-time LQS problems. By applying RL techniques, Jiang and Jiang \cite{Jiang2011} developed an ADP strategy to solve continuous-time optimal control problems where the systems subject to control-dependent noise.
~\\

However, to the author's best knowledge, there is no model-free results for continuous-time LQS optimal control problems where drift and diffusion terms contain both control and state variables.  The main contribution of this paper is that we propose a  model-free algorithm to solve this class of continuous-time LQS problems. 
~\\

To be specific, we propose a novel data-driven model-free PI algorithm to get the maximal solution to the SARE by using the data of the input and state collected on some time interval. The convergence proof of our model-free strategy is also been provided. 
~\\

The rest of the paper is organized as follows. In Section 2, the formulation of our problem and some preliminaries are presented. Section 3 develops our data-driven model-free PI algorithm. In Section 4, we provide a simulation example to illustrate the applicability of the proposed algorithm. In Section 5, some conclusions are presented.\\

\noindent{\bf Notation.} We denote the collections of non-negative integers, positive integers and real numbers by $\mathbb{Z}$, $\mathbb{Z^{+}}$ and $\mathbb{R}$. $\mathbb{R}^{n\times m}$ represents the collection of all $n\times m$ real matrices. $\mathbb{R}^{n}$ is the $n$-dimensional Euclidean space and $|\cdot|$ denotes its Euclidean norm for vector or matrix of proper size. Zero matrix (or vector) with appropriate dimension is denoted by $O$. We use $diag\{v\}$ to denote a square diagonal matrix whose main diagonal is the elements of vector $v$.  The sets of all symmetric matrices, positive definite matrices and semipositive definite matrices in $\mathbb{R}^{n\times n}$ are represented by $\textbf{S}^{n}$, $\textbf{S}^{n}_{++}$ and $\textbf{S}^{n}_{+}$, respectively. $w(\cdot)$ is a one-dimensional standard Brownian motion defined on a filtered probability space ($\Omega$,\,$\mathcal{F}$,\,$\{\mathcal{F}_{t}\}_{t\geqslant 0}$,\,$\mathbb{P}$) that satisfies usual conditions. 
Moreover, we use $\otimes$ to denote the Kronecker product and for any matrix $B\in\mathbb{R}^{m\times n}$, $vec(B)$ denotes a vectorization map from the matrix $B$ into a column vector of proper size, which stacks the columns of $B$ on top of one another, that is, $vec(B) = [b^T_1, b^T_2,\cdots, b^T_n]^T$, where $b_j\in\mathbb{R}^{n}$, $j=1,2,3,\cdots,n$, are the columns of $B$. For any $\xi \in\mathbb{R}^n$ and $F\in\mathbf{S}^n$, we define two operators as follows:
\begin{equation*}
\begin{split}
&vecs: \xi\in\mathbb{R}^l\rightarrow vecs(\xi)\in\mathbb{R}^{\frac{n(n+1)}{2}},\\
\text{and}\,\,\, 
&vech: F\in\mathbf{S}^l\rightarrow vech(F)\in\mathbb{R}^{\frac{n(n+1)}{2}},
\end{split}
\end{equation*}
where
\begin{small}
	\begin{equation*}
	\begin{split}
	&vecs(\xi)=[\xi_1^2, \xi_1\xi_2,\cdots,\xi_1\xi_n,x_2^2, x_2x_3,\cdots,\xi_{n-1}\xi_n,\xi_n^2]^T,\\ 
	&vech(F)\,=[f_{11}, 2f_{12},\cdots,2f_{1n},f_{22}, 2f_{23},\cdots,2f_{n-1,n},f_{nn}]^T,
	\end{split}
	\end{equation*}
\end{small}
and $\xi_j$, $j=1,2,\cdots,n$, is the $j$th element of $\xi$ and $f_{ji}$, $j,i=1,2,\cdots,n$, is the $(j,i)$th element of matrix $F$. For simplity, we denote $vecs(\xi)$ by $\overline{\xi}$ in this paper.\\

\section{PROBLEM FORMULATION}\label{sec2}
This section presents the formulation of our LQS optimal control problems. 
~\\

Consider a continuous-time time-invariant stochastic linear system as follows
\begin{equation}
\label{system}
\begin{cases}
\begin{split}
dx(s)= \,\,&[Ax(s)+Bu(s)]ds\\
&+[Cx(s)+Du(s)]dw(s),
\end{split}\\
x(0)=x_0,
\end{cases}
\end{equation}
where $x_0\in \mathbb{R}^{n}$ is the initial state. The cost functional is defined as
\begin{equation}\label{cost}
\begin{split}
J(u(\cdot))=\mathbb{E}\int_{0 }^{\infty}[x(s)^TQx(s)+u(s)^TRu(s)]ds,
\end{split}
\end{equation}
where $R>0$, $Q\geq0$ and  $[A,C|Q]$ is exactly detectable.

Now we give the definition of mean-square stabilizability.\\

\noindent{\bf Definition 1.} System (\ref{system}) is called mean-square stabilizable for any initial state $x_0$, if there exists a matrix $K\in \mathbb{R}^{m\times n} $ such that the solution of 
\begin{equation}
\label{eq4}
\begin{cases}
\begin{split}
dx(s)=\,\,&(A+BK)x(s)ds\\
&+(C+DK)x(s)dw(s), 
\end{split}\\
x(0)=x_0
\end{cases}
\end{equation}
satisfies $\lim _{s \to \infty}\mathbb{E}[x(s)^Tx(s)]=0$.
In this case, the feedback control $u(\cdot)=Kx(\cdot)$ is called stabilizing and the constant matrix $K$ is called a stabilizer of  system (\ref{system}). \\

\noindent{\bf Assumption 1.} System (\ref{system}) is mean-square stabilizable.\\

Under Assumption 1, we define the sets of admissible control as
\begin{equation}
\mathcal{U}_{ad}=\{u(\cdot)\in L^2_{\mathcal{F}}(\mathbb{R}^{m})|u(\cdot) \,\, \text{is\, stabilizing}\}.
\end{equation}

Our continuous-time LQS optimal control problems are given as follows:\\

\noindent{\bf Problem (LQS).} For any initial state $x_0\in \mathbb{R}^{n}$, we want to find an optimal control $u^*(\cdot) \in \mathcal{U}_{ad}$ such that
\begin{equation}
J(u^*(\cdot))=\inf \limits_{u(\cdot) \in \mathcal{U}_{ad}}J(u(\cdot)).
\end{equation}

Ni and Fang \cite{Ni2013} shows that the optimal control of Problem (LQS) can be obtained by solving the following stochastic algebraic Riccati equation (SARE)
\begin{equation}\label{SARE}
\begin{split}
PA&+A^TP+C^TPC+Q-(PB+C^TPD)\\
&\times(R+D^TPD)^{-1}(B^TP+D^TPC)=0.
\end{split}
\end{equation} 

Due to the nonlinear structure of SARE (\ref{SARE}), the analytical solution of (\ref{SARE}) is difficult to obtain. To our best knowledge, there are some iterative algorithms to get the approximation solution of (\ref{SARE}), one of which is the PI method developed in Ni and Fang \cite{Ni2013}. We summarize the method as the following lemma.\\

\noindent{\bf Lemma 1.} Assume $[A,C|Q]$ is exactly detectable. For a given stabilizer $K_0$,  let $P_{i} \in \textbf{S}^{n}_{+}$ be the solution of
\begin{equation}\label{eq5}
\begin{split}
&P_{i}(A+BK_i)+(A+BK_i)^TP_{i}+Q\\
&+(C+DK_i)^TP_{i}(C+DK_i)+K_i^TRK_i=0,
\end{split}
\end{equation}
where $K_{i}$ is updated by
\begin{equation}\label{eq6}
\begin{split}
K_{i+1}=-(R+D^TP_{i}D)^{-1}
(B^TP_{i}+D^TP_{i}C).
\end{split}
\end{equation}

Then $P_i$ and $K_i$, $i=0,1,2,3,\cdots$ can be uniquely determined at each iteration step, and the following conclusions hold:\\
(i) $K_i$, $i=0,1,2,\cdots$, are stabilizers.\\
(ii) $\lim _{i \to \infty}P_i=P^*$, $\lim _{i \to \infty}K_i=K^*$, where $P^*$ is a nonnegative definite solution to SARE (\ref{SARE}) and $K^*=-(R+D^TP^*D)^{-1}(B^TP^*+D^TP^*C)$. In this case,  $ u^*(\cdot)=K^*x^*(\cdot)$ is an optimal control of Problem (LQS).\\

Note that the above method needs all knowledge of the system matrices, which are difficult to obtain in the real world. Thus we want to develop a model-free algorithm to solve $P_{i}$ and $K_{i}$ without using the information of the coefficient matrices $A$, $B$, $C$, $D$ in system (\ref{system}).\\

\section{MODEL-FREE PI ALGORITHM }\label{sec3}

In this section, we present our data-driven PI algorithm that does not rely on all knowledge of the coefficient matrices in system (\ref{system}).\\ 

To this end, we first rewrite (\ref{eq5}) as

\begin{equation}\label{eq7}
\begin{split}
&A^TP_{i}+P_{i}A+C^TP_{i}C\\
=&-P_iBK_i-K_i^TB^TP_i-Q-K_i^TD^TP_iC\\
&-C^TP_iDK_i-K_i^TD^TP_iDK_i-K_i^TRK_i.\\
\end{split}
\end{equation}

Then, by Ito's formula, we know 
\begin{equation}\label{eq8}
\begin{split}
&d\big(x(s)^TP_{i}x(s)\big)\\
=\bigg\{&x(s)^T\big(A^TP_{i}+P_{i}A+C^TP_{i}C\big)x(s)\\
&+2u(s)^T\big(B^TP_{i}+D^TP_{i}C\big)x(s)\\
&+u(s)^TD^TP_{i}Du(s)\bigg\}ds+\bigg\{\cdots\bigg\}dw(s).
\end{split}
\end{equation}\\

Combining it with (\ref{eq7}), we have
\begin{equation}\label{eq}
\begin{split}
&d\big(x(s)^TP_{i}x(s)\big)\\
=\bigg\{&-x(s)^T\big(Q+K_i^TRK_i\big)x(s)\\
&+2\big(u(s)-K_ix(s)\big)^T\big(B^TP_{i}+D^TP_{i}C\big)x(s)\\
&+u(s)^TD^TP_{i}Du(s)\\
&-x(s)^TK_i^TD^TP_{i}DK_ix(s)\bigg\}ds+\bigg\{\cdots\bigg\}dw(s).\\
\end{split}
\end{equation}
Integrating (\ref{eq}) from $t$ to $t+\triangle  t$ and taking expection $\mathbb{E}$, we get 

\begin{equation}\label{eq10}
\begin{split}
&\mathbb{E}\big[x(t+\triangle  t)^TP_{i}x(t+\triangle  t)-x(t)^TP_{i}x(t)\big]\\
&-2\mathbb{E}\int_{t}^{t+\triangle  t}\big(u(s)-K_ix(s)\big)^TM_{i}x(s)ds\\	
&-\mathbb{E}\int_{t}^{t+\triangle  t}u(s)^TH_{i}u(s)ds\\
&+\mathbb{E}\int_{t}^{t+\triangle  t}x(s)^TK_i^TH_{i}K_ix(s)ds\\
=&-\mathbb{E}\int_{t}^{t+\triangle  t}x(s)^T\big(Q+K_i^TRK_i\big)x(s)ds,		
\end{split}
\end{equation}
where $M_{i}=B^TP_{i}+D^TP_{i}C$, $H_{i}=D^TP_{i}D$, $t\geq 0$, $\triangle t$ is any positive real number and $x(\cdot)$ is governed by system (\ref{system}) with any control $u(\cdot)$.\\

Next, we give some symbols to develop our data-driven model-free PI algorithm. We define matrices $\eta  _{\overline{x}} \in \mathbb{R}^{q\times \frac{n(n+1)}{2}} $, $\eta  _{\overline{u}}\in \mathbb{R}^{q\times \frac{m(m+1)}{2}} $, $\eta  _{\overline{K_ix}}\in \mathbb{R}^{q\times \frac{m(m+1)}{2}}$, $i=0,1,2,\cdots$, $\eta  _{xu} \in \mathbb{R}^{q\times mn} $ and $\eta  _{xx} \in \mathbb{R}^{q\times n^2} $,  as follows
\begin{equation*}
\begin{split}
\eta  _{\overline{x}}=\mathbb{E}\bigg[\overline{x(t_1)} -\overline{x(t_0)},\cdots,
\overline{x(t_q)}-\overline{x(t_{q-1})}\bigg]^T,\\
\end{split}
\end{equation*}
\begin{equation*}
\begin{split}
\eta  _{\overline{u}}=\mathbb{E}\bigg[\int_{t_0}^{t_1} \overline{u(s)}ds,\cdots,
\int_{t_{q-1}}^{t_q}\overline{u(s)}ds\bigg]^T,\\
\end{split}
\end{equation*}
\begin{equation*}
\begin{split}
\eta  _{\overline{K_ix}}=\mathbb{E}\bigg[\int_{t_0}^{t_1} \overline{K_ix(s)}ds,\cdots,
\int_{t_{q-1}}^{t_q}\overline{K_ix(s)}ds\bigg]^T,\\
\end{split}
\end{equation*}
\begin{equation*}
\begin{split}
\eta  _{xx}=\mathbb{E}\bigg[\int_{t_0}^{t_1}x(s)\otimes x(s)ds,\cdots,
\int_{t_{q-1}}^{t_q}x(s)\otimes x(s)ds\bigg]^T,\\
\end{split}
\end{equation*}
\begin{equation*}
\begin{split}
\eta  _{xu}=\mathbb{E}\bigg[\int_{t_0}^{t_1}x(s)\otimes u(s)ds,\cdots,
\int_{t_{q-1}}^{t_q}x(s)\otimes u(s)ds\bigg]^T,\\
\end{split}
\end{equation*}
where $q\in\mathbb{Z^{+}}$ is any positive integer and $0\leq t_0<t_1<t_2<\cdots<t_q$.\\

For any given $K_i$, (\ref{eq10}) implies 
\begin{equation}\label{solve1}
\Psi_i\begin{bmatrix}
vech(P_{i})\\
vec(M_{i})\\
vech(H_{i})\\
\end{bmatrix}=\Theta_i, 
\end{equation}
where $\Psi_i \in \mathbb{R}^{q\times (\frac{n(n+1)}{2}+mn+\frac{m(m+1)}{2})}$ and $\Theta_i \in \mathbb{R}^{q}$ are defined as
\begin{equation*}
\begin{split}
\Theta_i=\big[-\eta _{xx}vec(Q+K_i^TRK_i)\big],
\end{split}
\end{equation*}
\begin{equation*}
\begin{split}
\Psi_i=\big[\eta  _{\overline{x}},2\eta _{xx}(I_n\otimes K_i^T)-2\eta _{xu}, \eta  _{\overline{K_ix}}-\eta  _{\overline{u}}  \big].
\end{split}
\end{equation*}

If $\Psi_i$ has full column rank for any $i\in\mathbb{Z}$,  (\ref{solve1}) can be directly transformed to 
\begin{equation}\label{solve2}
\begin{bmatrix}
vech(P_{i})\\
vec(M_{i})\\
vech(H_{i})\\
\end{bmatrix}=(\Psi_i^T\Psi_i)^{-1}\Psi_i^T\Theta_i.
\end{equation}\\

Next, we show that, under condition (\ref{rank}) in the following lemma, $\Psi_i$, $i=0,1,\cdots$, has full column rank.\\

\noindent{\bf Lemma 3.} If there exists a $q_0 \in \mathbb{Z^{+}}$, such that, for all $q \geq q_0$,
\begin{equation}\label{rank}
rank([\eta _{xx},\,\,\eta _{xu},\,\,\eta  _{\overline{u}} ])=\frac{n(n+1)}{2}+mn+\frac{m(m+1)}{2},
\end{equation}
then, $\Psi_i$, $i=0,1,\cdots$, has full column rank.\\

\noindent{\bf Proof.} It is enough to prove that
\begin{equation}\label{eq14}
\Psi_iV=O,\,\,\,\forall i\in\mathbb{Z},
\end{equation} 
has the unique solution $V=O$, where $O$ is a zero matrix (or vector) with appropriate dimension and  $V\in\mathbb{R}^{mn+\frac{n(n+1)}{2}+\frac{m(m+1)}{2}}$ .\\

To achieve it, we now prove it by contradiction. We assume $V=[vech(N)^T,vec(F)^T,vech(G)^T]^T\in  \mathbb{R}^{mn+\frac{n(n+1)}{2}+\frac{m(m+1)}{2}}$ is a nonzero column vector, where $vech(N)\in \mathbb{R}^{\frac{n(n+1)}{2}}$,  $vec(F)\in \mathbb{R}^{mn}$ and $vech(G)\in\mathbb{R}^{\frac{m(m+1)}{2}}$. Then, by the definitions of $vech(\cdot)$ and $vec(\cdot)$, two symmetric matrices $N\in\textbf{S}^{n}$, $G\in\textbf{S}^{m}$ and a matrix $F\in\mathbb{R}^{m\times n}$ can be uniquely determined by $vech(N)$, $vech(G)$ and $vec(F)$, respectively.\\

Applying Ito's formula to $x(s)^TNx(s)$, we derive 
\begin{equation}\label{eq13}
\begin{split}
&\mathbb{E}\big[x(t+\triangle  t)^TNx(t+\triangle  t)-X(t)^TNX(t)\big]\\
=
\,\,&\mathbb{E}\int_{t}^{t+\triangle  t}x(s)^T\big(A^TN+NA+C^TNC\big)x(s)ds\\
&+2\mathbb{E}\int_{t}^{t+\triangle  t}u(s)^TB^TNx(s)ds\\
&+2\mathbb{E}\int_{t}^{t+\triangle  t}u(s)^TD^TNCx(s)ds\\
&+\mathbb{E}\int_{t}^{t+\triangle  t}u(s)^TD^TNDu(s)\big)ds,\\		
\end{split}
\end{equation}
where $x(\cdot)$ is governed by system (\ref{system}) with the same input $u(\cdot)$ as in (\ref{eq10}).\\

Using (\ref{eq10}), (\ref{eq13}) and the definition of $\Psi_i$, we have
\begin{equation}\label{eq15}
\Psi_iV=\eta _{xx}vec(\mathcal{T})+\eta _{xu}vec(\mathcal{J})+\eta _{\overline{u}}vech(\mathcal{L}),
\end{equation}
where
\begin{equation}\label{eq16}
\begin{split}
\mathcal{T}=&A^TN+NA+C^TNC+K_i^TGK_i\\
&+K_i^TF+F^TK_i\\		
\end{split}
\end{equation}
\begin{equation}\label{eq17}
\begin{split}
\mathcal{J}=2B^TN+2D^TNC-2F,
\end{split}
\end{equation}
\begin{equation}\label{eq18}
\begin{split}
\mathcal{L}=D^TND-G.
\end{split}
\end{equation}

Since $\mathcal{T}$ is a symmetric matrix, we get
\begin{equation}
\eta _{xx}vec(\mathcal{T})=I_{\overline{x}}vech(\mathcal{T}),
\end{equation}
where $I_{\overline{x}}\in\mathbb{R}^{q\times\frac{n(n+1)}{2}}$ and
\begin{equation}
\begin{split}
I_{\overline{x}}=\mathbb{E}\bigg[\int_{t_0}^{t_1}\overline{x(s)}ds,\cdots,
\int_{t_{q-1}}^{t_q}\overline{x(s)}ds\bigg]^T.\\
\end{split}
\end{equation}

Then, (\ref{eq14}) and (\ref{eq15}) imply 
\begin{equation}\label{eq19}
[I_{\overline{x}},\eta _{xu},\eta _{\overline{u}}]\begin{pmatrix}
vech(\mathcal{T})\\
vec(\mathcal{J})\\
vech(\mathcal{L})\\
\end{pmatrix}=O.
\end{equation}

It is easy to see that $[I_{\overline{x}},\eta _{xu},\eta _{\overline{u}}]$ has full column rank under condition (\ref{rank}). Then, the solution to (\ref{eq19}) is $vech(\mathcal{T})=O$, $vec(\mathcal{J})=O$ and $vech(\mathcal{L})=O$, and thus $\mathcal{T}=O,\mathcal{J}=O$ and $\mathcal{L}=O$.\\

Next, since $K_i$ is a stabilizer, by Definition 1, we know the trajectory of 
\begin{equation}
\label{system2}
\begin{cases}
\begin{split}
dx(s)= \,\,&\Big[(A+BK_i)x(s)\Big]ds\\
&+\Big[(C+DK_i)x(s)\Big]dw(s),\\
\end{split}\\
x(0)=x_0\in\mathbb{R}^n
\end{cases}
\end{equation}
satisfies $\lim_{s\rightarrow+\infty}\mathbb{E}\big[x(s)^Tx(s)\big]=0$.\\

For any $t>0$, applying Ito's formula to $d\big(x(s)^TNx(s)\big)$, we get 
\begin{equation}\label{eq88}
\begin{split}
&\mathbb{E}\big[x^T(t)Nx(t)\big]-x_0^TNx_0\\
=\,\,&\mathbb{E}\int_{0}^{t}x^T(s)\big((A+BK_i)^TN+N(A+BK_i)\\
&+(C+DK_i)^TN(C+DK_i)\big)x(s)ds,\\	
\end{split}
\end{equation}
where $x(\cdot)$ is governed by (\ref{system2}).\\

Then, by (\ref{eq16}), (\ref{eq17}), (\ref{eq18}), $\mathcal{T}=0$, $\mathcal{J}=0$ and $\mathcal{L}=0$, we can easily see from (\ref{eq88}) that
$\mathbb{E}\big[x^T(t)Nx(t)\big]-x_0^TNx_0=0$. Letting $t\rightarrow+\infty$, we have $x_0^TNx_0=\lim_{t\rightarrow+\infty}\mathbb{E}\big[x^T(t)Nx(t)\big]=0$. Notice that $x_0$ can be any element in $\mathbb{R}^n$ and $N\in\mathcal{S}^n$, we have $N=0$. Then it follows from (\ref{eq16}), (\ref{eq17}), (\ref{eq18}), $\mathcal{T}=0$, $\mathcal{J}=0$ and $\mathcal{L}=0$ that $G=0$ and $F=0$, which contradicts with $V\neq0$. The proof is completed.$\hfill\blacksquare$\\

Using above notations, our model-free algorithm is given in Algorithm 1.\\

\begin{algorithm}[h]
	\caption{}
	\label{A2}
	\begin{algorithmic}[1]
		
		\State Initial $i=0$ and select $K_0$ as a stabilizer for system (\ref{system}). Take $u(\cdot)=K_0x(\cdot)+e(\cdot)$ as the input to system (\ref{system}) on time interval $[t_0,t_q]$, where $e(\cdot)$ is the exploration noise. Calculate $\eta  _{\overline{x}}$, $\eta  _{\overline{u}} $, $\eta  _{xu} $ and $\eta  _{xx}$.
		
		\State \textbf{repeat}
		
		\State  Compute $\eta  _{\overline{K_ix}}$ and solve $P_{i}$, $M_{i}$ and $H_{i}$ from (\ref{solve2}).

		\State  $K_{i+1}=-(R+H_{i})^{-1}M_{i}$.
		
		\State $i\leftarrow i+1$.
		\State \textbf{Until} {$|P_{i+1}-P_{i}|<\varepsilon$.}
	\end{algorithmic}
\end{algorithm}

Finally, we show the convergence of our algorithm.\\

\noindent{\bf Theorem 1.} Under rank condition (\ref{rank}), starting from a stabilizer $K_0$, the sequences $\{P_i\}_{i=0}^{\infty}$ and $\{K_i\}_{i=1}^{\infty}$ obtained from Algorithm 1 satisfy $\lim_{i\rightarrow \infty}P_i=P^*$ and $\lim_{i\rightarrow \infty}K_i=K^*$.\\

\noindent{\bf Proof.} Given a stabilizer $K_i$, if $P_{i}\in\textbf{S}^{n\times n}$ is the solution of (\ref{eq5}), $M_{i}$ and $H_{i}$ can be uniquely determined by $M_{i}=B^TP_{i}+D^TP_{i}C$ and $H_{i}=D^TP_{i}D$, respectively. Thus, (\ref{eq10}) implies that $P_{i}$, $M_{i}$ and $H_{i}$ must satisfy (\ref{solve2}).\\

Moreover, if (\ref{rank}) holds, (\ref{solve2}) has the unique solution $(P_{i},M_{i},H_{i})$. Otherwise, (\ref{solve2}) has two different solutions and thus contradicts with rank condition (\ref{rank}). \\

Therefore, under condition (\ref{rank}), $P_{i}$ and $K_{i}$, $i=0,1,2,\cdots$, obtained from Algorithm 1 are equivalent to the solution of (\ref{eq5}) and (\ref{eq6}). Then the convergence of the proposed algorithm can be guaranteed by Lemma 1. $\hfill\blacksquare$\\

\section{NUMERICAL EXAMPLE}\label{sec4}

This section will present a simulation example to illustrate the  feasibility of Algorithm 1.\\

We consider system (\ref{system}) with $n=2$ and $m=1$, 
\begin{equation*}
A=
\begin{bmatrix}

0  & -0.6\\
0.6   & -0.3
\end{bmatrix},
B=
\begin{bmatrix}

0.05\\
0.01
\end{bmatrix},
\end{equation*} 
\begin{equation*}
C=
\begin{bmatrix}
-0.02 &   0.03\\
-0.05    &0.02
\end{bmatrix},
D=
\begin{bmatrix}
0.001\\
0.03    
\end{bmatrix},
\end{equation*} 
and $x_0=[0.5,-0.1]^T$.  The weighting matrices in the cost functional are choosed as $R=1>0$ and $Q=diag(1,0.5)\geq 0$.\\  

By implementing Algorithm 1, we can obtain
\begin{equation*}
\widetilde{P}^*=
\begin{bmatrix}
2.9072352 & -0.8296538\\
-0.8296538&2.4975686
\end{bmatrix},
\end{equation*}
\begin{equation*}
\widetilde{K}^*=
\begin{bmatrix}
-0.0669434 & 0.0064058\\
\end{bmatrix}.
\end{equation*}

Moreover, to check the error of the proposed algorithm, we denote the left sides of (\ref{SARE}) and (\ref{eq5}) as $\mathcal{R}_1(P)$ and $\mathcal{R}_2(P, K)$. Then we have $|\mathcal{R}_1(\widetilde{P}^*)|=2.0820041\times10^{-3}$ and $|\mathcal{R}_2(\widetilde{P}^*,\widetilde{K}^*)|=2.0833488\times10^{-3}$.\\

\section{CONCLUSION}\label{sec5}

This paper has developed a model-free PI algorithm to solve infinite-horizon LQS problems, i.e., Problem (LQS). By applying ADP techniques, the solution of Problem (LQS) can be learned from the collected data. Moreover, an example is given to show the applicability of the obtained algorithm. \\

\section*{Acknowledgment}
\noindent The author is grateful for the constructive comments of  Professor Guangchen Wang, which leads to an improvement of this work.

\end{document}